\documentclass{amsart}
\theoremstyle{definition}

\theoremstyle{remark}

\newcommand{\R}{\mathbb R}

\newcommand{\Z}{\mathbb Z}
\renewcommand\Re[1]{\operatorname{Re}}
\newcommand{\norm}[1]{\left\Vert#1\right\Vert}
\newcommand{\biggnorm}[1]{\biggl\Vert#1\biggr\Vert}
\def\abs|#1|{\left\vert#1\right\vert}

\newcommand\mcD{\mathcal{D}}
\newcommand\mcI{\mathcal{I}}
\newcommand\mcJ{\mathcal{J}}
\newcommand\wrt{\,{\rm d}\/}

\begin{document}

\title{Some light on Littlewood--Paley theory}
\author{Michael Cowling}
\address{School of Mathematics, University of New South Wales,
Sydney NSW 2052, AUSTRALIA}
\email{m.cowling@unsw.edu.au}

\author{Terence Tao}
\address{Department of Mathematics, University of California at Los Angeles,
Los Angeles CA 90095, USA}
\email{tao@math.ucla.edu}
\thanks{}
\date{}
\dedicatory{}
\commby{}


\begin{abstract}
The purpose of this note is to correct an error in a paper of
M.~Cowling, G. Fendler and J.J.F. Fournier, and to give a
counterexample to a conjecture of J.-L. Rubio de Francia.
\end{abstract}

\maketitle
Classical Littlewood--Paley theory (in the context of analysis on
$\R$) deals with expressions of the form
\[
\Bigl( \sum_{I\in \mcI} \abs|S_If|^2\Bigr)^{1/2},
\]
where $\mcI$ is a collection of (essentially) disjoint intervals in
$\R$, typically the set $\mcD$ of dyadic intervals $\{ [-2^{n+1},
-2^{n}], [2^{n},2^{n+1}]: n \in \Z\}$, and $S_I$ is the operator of Fourier
multiplication by the characteristic function of the interval $I$,
i.e., $(S_I f)\hat{\phantom{f}} = \chi_I \hat f$.  The classical
Littlewood--Paley inequality states that, if $1<p<\infty$, then
there exist (positive) constants~$A_p$ and $B_p$ such that
\begin{equation*}\label{f: LP}
     A_p \norm{f}_p
\leq     \biggnorm{\Bigl( \sum_{I\in \mcD} \abs|S_If|^2\Bigr)^{1/2}}_p
\leq B_p \norm{f}_p
\qquad\forall f \in L^p(\R).
\end{equation*}
The collection of intervals may be replaced by other collections of
intervals, but for the left hand inequality to hold, it is
important that the union of the intervals be (essentially) all
of~$\R$; this is not a restriction for the right hand inequality.

This has been generalized in a number of ways.  J.-L.~Rubio de
Francia \cite{RdF} proved that, if $\mcI$ is any collection of
disjoint intervals and $2 \leq p < \infty$, then an inequality
\begin{equation}\label{f: RdF LP}
\biggnorm{\Bigl( \sum_{I\in \mcD} \abs|S_If|^2\Bigr)^{1/2}}_p
\leq B_p \norm{f}_p
\qquad\forall f \in L^p(\R)
\end{equation}
still holds.  He also observed that this inequality cannot hold if
$1<p<2$.  Indeed, if $\mcJ$ is the collection of all intervals
$\{[n,n+1]: n\in\Z\}$, and $f_N$ denotes the function on $\R$ whose
Fourier transform~$\hat f_N$ is the characteristic function
$\chi_{[0,N]}$, for some positive integer~$N$, then it is
straightforward to check that
\[
\biggnorm{\Bigl( \sum_{I\in \mcJ} \abs|S_I f_N|^2\Bigr)^{1/2}}_p
= N^{1/2} \norm{f_1}_p,
\]
while
\[
\norm{f_N}_p = N^{1/p'} \norm{f_1}_p,
\]
where $p'$ is the dual index to $p$, that is, $p' = p/(p-1)$. Thus
\eqref{f: RdF LP} can hold only if $p\geq 2$.  However, for this
example, the modified Littlewood--Paley inequality
\begin{equation}\label{f: new RdF LP}
\biggnorm{\Bigl( \sum_{I\in \mcD} \abs|S_If|^{p'}\Bigr)^{1/p'}}_p
\leq B_p \norm{f}_p
\end{equation}
holds.  Perhaps a little optimistically, Rubio de Francia
\cite[p.~10]{RdF} conjectured that \eqref{f: new RdF LP} might
always hold.  One of the aims of this paper is to provide a
counterexample to this conjecture.

Shortly after Rubio de Francia's paper, Cowling, Fendler and
Fournier \cite{CFF} investigated variants of Littlewood--Paley
theory in which mixed norms like that in \eqref{f: new RdF LP}
appear.  These were used to give some examples of multipliers with
some special properties.  Cowling, Fendler and Fournier
\cite[p.~340]{CFF} used the space called $D(\R)$ of all
integrable functions~$f$ on~$\R$ such that $\int_n^{n+1} f(x) \wrt
x = 0$ for all integers~$n$.  In particular, they claimed 
that the real interpolation space $[D(\R), L^2(\R)]_{\theta,p}$ 
is the Lebesgue space~$L^p(\R)$, where $1/p = 1 -\theta/2$.  
A second purpose of this paper is to disprove this
assertion; consequently all results based on this ``fact'' are
suspect.

This paper owes much to T.-S.~Quek, who observed that the
interpolation theorem above would (if true) imply 
Rubio de Francia's conjecture.

\subsection{A counterexample to Rubio de Francia's conjecture}

Let $\mcI$ be the family of all intervals~$I_{j,n}$ of the form
\[
[n + j 2^{-n} , n + (j+1) 2^{-n}],
\]
where $n = 0, 1, 2, \dots$ and $0 \leq j < 2^n$. Again take the
function~$f_N$ such that $\hat f = \chi_{[0,N]}$; then
\begin{equation}\label{f: one side}
\norm{f_N}_p = N^{1/p'}\norm{f_1}_p.
\end{equation}

Now consider one of the intervals $I_{j,n}$ above, where $n < N$.
The absolute value of the function~$S_{I_{j,n}} f_N$ is equal to
the absolute value of the function~$f_{2^{-n}}$ whose Fourier
transform is the characteristic function~$\chi_{[0, 2^{-n}]}$,
i.e., to the absolute value of the function $x \mapsto\sin(2^{-n}
\pi x) /\pi x$ (using the Fourier transform with $2\pi$ in the
exponent). Thus $\abs|S_{I_{j,n}} f_N|$ is greater than
$2^{1-n}/\pi$ on the interval of length $2^n$ centered at the
origin. But for each $n$ there are $2^n$ intervals $I_{j,n}$.
Summing, we see that 
\[
\Bigl(\sum_{j=0}^{2^n-1} \abs|S_{I_{j,n}} f_N(x)|^{p'} \Bigr)^{1/p'}
\geq \frac{2^{1-n/p}}{\pi} \chi_{[-2^{n-1},2^{n-1}]}(x)
\]
for all $x$ in $\R$, whence
\[
\Bigl(\sum_{I\in\mcI} \abs|S_{I} f_N(x)|^{p'}\Bigr)^{1/p'}
\geq \frac{2}{\pi (4\abs|x|)^{1/p}}
\]
when $1/4 \leq \abs|x| \leq 2^{N}/4$, and so
\[
\biggnorm{ \Bigl(\sum_{I\in\mcI} |S_I f|^{p'}\Bigr)^{1/p'} }_p
\geq \frac{2}{\pi} \Bigl( 2\int_{1/4}^{2^{N}/4} 
                        \frac{1}{4x} \wrt x \Bigr)^{1/p}
 =   \frac{2}{\pi} \Bigl( \frac{N \log 2}{2} \Bigr)^{1/p}.
\]
This inequality, together with \eqref{f: one side}, shows that
Rubio's conjecture cannot hold if $1<p<2$.

\subsection{The reason why the interpolation ``theorem'' is not correct}
A function $u$ lies in the real interpolation space 
$[D(\R), L^2(\R)]_{\theta,q}$ (constructed using the $J$-method)
if and only if it may be represented as a vector-valued integral, 
with values in $D(\R)\cap L^2(\R)$, convergent in $D(\R) + L^2(\R)$,
and \emph{a fortiori} convergent in $L^1(\R)+L^2(\R)$:
\[
u = \int_0^\infty u_t\,\frac{\wrt t}{t},
\]
where
\[
\biggl( \int_0^\infty \bigl( 
        t^{-\theta} \max\{\norm{u_t}_{D(\R)}, t\norm{u_t}_{L^2(\R)} \} 
                        \bigr)^q \, \frac{\wrt t}{t}\biggr)^{1/q} 
        < \infty.
\]
For almost all $t$ in $\R^+$, $u_t$ has to lie in $D(\R)$, so
\[
\int_n^{n+1} u_t (x) \wrt x = 0.
\]
Since the map $M_n: u \mapsto \int_n^{n+1} u(x) \wrt x$ is 
continuous on $L^1(\R)+L^2(\R)$, 
\[
M_n u = \int_0^\infty M_n u_t\,\frac{\wrt t}{t} = 0,
\]
so the interpolation space is a subspace of the space of all 
locally integrable functions on $\R$ whose integrals on  the
intervals $[n,n+1]$ all vanish.



\begin{thebibliography}{9}
\bibitem{CFF} M. Cowling, G. Fendler and J.J.F. Fournier, 
Variants on Littlewood--Paley theory,
\emph{Math.\ Annalen}, \textbf{285} (1989), 333--342.

\bibitem{RdF} J.-L. Rubio de Francia, 
A Littlewood--Paley theorem for arbitrary intervals,
\emph{Rev.\ Mat. Iberoamericana}, \textbf{1} (1985), 1--14.

\end{thebibliography}
\end{document}